\documentclass[11pt]{amsart}
\usepackage{amsmath}
\usepackage{amsthm}
\usepackage{graphicx} 
\usepackage{rotate} 
\usepackage{pxfonts}
\usepackage{txfonts}
\usepackage{url}

\theoremstyle{plain}

\theoremstyle{plain}
\newtheorem{theorem}{Theorem}[section]
\newtheorem{proposition}[theorem]{Proposition}
\newtheorem{lemma}[theorem]{Lemma}

\theoremstyle{definition}

\newtheorem{question}[theorem]{Question}

\newcommand{\Fs}{\mathcal F^s}
\newcommand{\Fu}{\mathcal F^u}

\renewcommand{\int}[1]{\mathring{#1}}

\newcommand{\N}{\mathbb{N}}

\newcommand{\R}{\mathbb{R}}

\newcommand{\Sph}{\mathbb{S}}
\newcommand{\SLZ}{\mathrm{SL}_2(\Z)}

\newcommand{\Z}{\mathbb{Z}}

\begin{document}

\title[Small dilatation and Lorenz knots]{Small dilatation homeomorphisms \\as monodromies of Lorenz knots}

\author{Pierre Dehornoy}
\address{Institut Fourier, Universit\'e Joseph Fourier-Grenoble 1, 100 rue des math\'ematiques, BP 74, 38402 Saint-Martin-d'H\`eres, France}
\email{pierre.dehornoy@ujf-grenoble.fr}
\thanks{Thanks to Eko Hironaka and Ruth Kellerhals for organizing the conference \emph{Growth and Mahler measures in geometry and topology} in July 2013 and to the Institut Mittag-Leffler for hosting this conference. \\Partially supported by SNF project 137548 \emph{Knots and surfaces}.}

\date{\today}

\begin{abstract}
We exhibit low-dilatation families of surface homeomorphisms among monodromies of Lorenz knots. 
\end{abstract}


\maketitle

Pseudo-Anosov homeomorphisms are topological/dynamical objects that can be seen as geometric counterparts of non-cyclotomic irreducible polynomials. 
In this dictionnary, Mahler measure becomes what is called \emph{geometrical dilatation}. 
A natural task is then to exhibit (or better, to classify) homeomorphisms with low dilatation. 
There exist several constructions of such low-dilatation families (see the census~\cite{Hironaka11}): for example using fibered faces of the Thurston norm ball~\cite{McMullenTeich}, or using mixed-sign Coxeter diagrams~\cite{Hironaka12}.
The goal of this note is to exhibit an additional construction, that comes from \emph{Lorenz knots}, that is, periodic orbits of the Lorenz vector field.

This text contains few new results. 
Most of the content comes from the article~\cite{DLorenz} where we investigated the homological dilatation of Lorenz knots. 
The interest here is $(i)$ to restrict our attention to subfamilies of Lorenz knots for which a stronger statement can be obtained with much less technical efforts, $(ii)$ to notice that what was proven for homological dilatation in~\cite{DLorenz} can also be proven for geometrical dilatation.


\section{Introduction}

\begin{figure}[!ht]
	\includegraphics*[width=\textwidth]{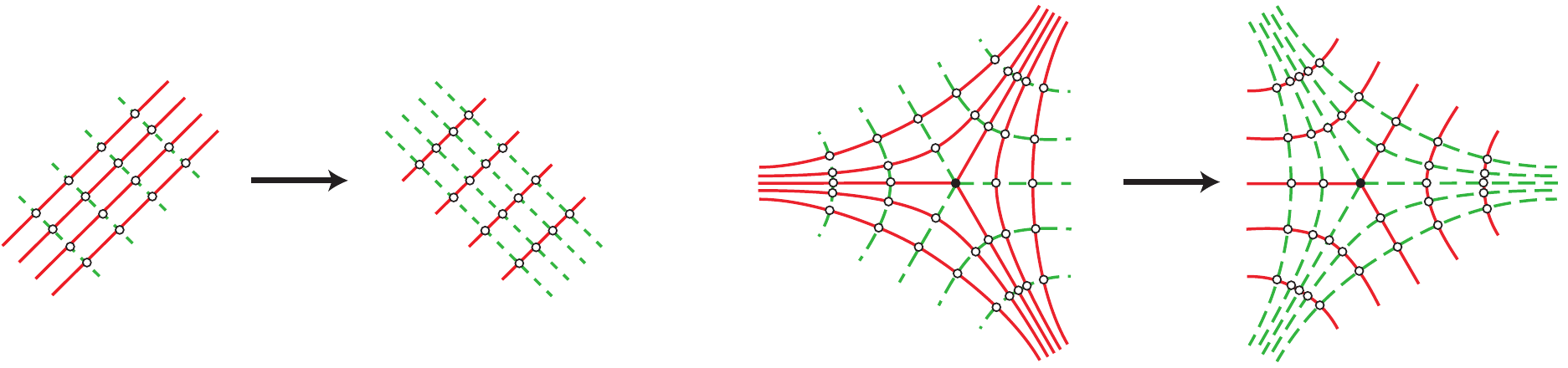}
	\includegraphics*[width=\textwidth]{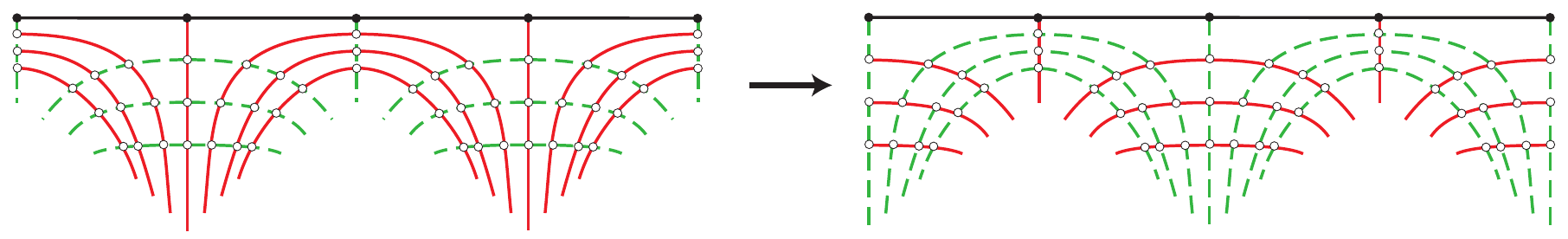}
	\caption{\small The stable and unstable foliations of a pseudo-Anosov homeomorphism (in red and green respectively). On the top left around a regular point, on the top right around a 3-prong singularity, and on the bottom around a boundary component with several singularities.}
  	\label{F:Singularities}
\end{figure}

It is known since Thurston~\cite{FLP, ThurstonSurfaces} that every homeomorphism of a surface is isotopic to either a periodic homeomorphism, or to a pseudo-Anosov one, or to a reducible one.
A \emph{pseudo-Anosov homeomorphism} of a surface~$S$ is a homeomorphism~$h$ such that $S$ admits two transverse measured foliations, called \emph{stable} and \emph{unstable} and usually denoted by~$(\Fs, \mu^s),$ and $(\Fu, \mu^u)$, that are invariant under~$h$, and such that there exists a positive real~$\lambda(h)$, called the \emph{geometrical dilatation} of~$h$, such that $\mu^s$ and $\mu^u$ are uniformly multiplied by $\lambda(h)^{-1}$ and $\lambda(h)$ under~$h$ respectively. 
Another property of~$h$ is that all closed curves on~$S$ are stretched at speed~$\lambda(h)$: for any auxiliary metric on~$S$ and for any closed curve~$\gamma$ on~$S$, we have~$\lim_{n\to\infty}\log(\Vert h^n(\gamma)\Vert)/n=\log(\lambda(h))$.
The most standard example is given by the action of a hyperbolic matrix in~$\SLZ$ on the torus~$\R^2/\Z^2$. In this case the invariant foliations are given by the two eigendirections of the matrix and the dilatation is the largest eigenvalue. 
On surfaces of higher genus, the foliations have prong-type singularities (see Figure~\ref{F:Singularities}).
A reducible homeomorphism is one that admits invariant curves. 
Such curves divide the surface into elementary pieces where the dynamics is either periodic or of pseudo-Anosov type.

This decomposition result of Thurston can be compared with the fact that every polynomial is either cycloctomic, or has positive Mahler measure, or is reducible.
The dilatation factor for a pseudo-Anosov homeomorphism, or better its normalized version $\lambda(h)^{\vert\chi(S)\vert}$, is a natural counterpart of the Mahler measure.
In particular, on a given surface, it is easy to find homeomorphisms with arbitrarily high dilatation (for example by iterating a fixed homeomorphism), but low-dilatation homeomorphisms are harder to construct. 
For $g\ge 1$ and $s\ge 0$, one usually defines $\delta_{g,s}$ as the infimum of the dilatation of a pseudo-Anosov homeomorphism on a closed surface of genus~$g$ with $s$ boundary components. 
It is \emph{a priori} not clear whether $\delta_{g,s}$ is 1 or larger, and in the latter case whether it is a minimum or not.
D.\,Penner showed~\cite{Penner} that $\delta_{g,s}$ is actually a minimum larger than 1, and that there exist two positive constansts $c_1, c_2$ such that for all $g$, one has $\frac{c_1}g\le \log \delta_{g,0}\le \frac{c_2}g$ (similar results hold for $s\ge 1$).
The optimal value of $c_1$ is not known, and the only known values of $\delta_{g,0}$ are $\delta_{1,0}$ and~$\delta_{2,0}$~\cite{CH}. 
It also follows from the work of Penner that, for any positive $D$ and for a fixed surface, only finitely many mapping classes have a dilatation smaller than $D$. 
It is then natural to study what happens when the genus tends to infinity.

A family $(h_n, S\!_n)_{n\in\N}$, where $S\!_n$ is a closed orientable surface and $h_n$ a pseudo-Anosov homeomorphism of~$S\!_n$, is said to be of \emph{low-dilatation} if the sequence $\log (\lambda(h_n))\vert\chi(S\!_n)\vert$ is bounded. 

Low-dilatation families are well understood in the context of 3-manifolds, where a theorem of B.\,Farb, C.\,Leininger and D.\,Margalit~\cite{FLM} states that the punctured suspensions of a low-dilatation family live in some fibered faces of the Thurston norm ball of a \emph{finite} number of 3-manifolds.
However it is still unknown 
how different homeomorphisms having the same suspensions are related. 

\begin{question}[Farb-Leininger-Margalit~\cite{FLM}]
\label{Q:QPeriodic}
Given a positive number~$D$, what can be said about the dynamics of those homeomorphisms with normalized dilatation smaller than~$D$ ({\it i.e.}, those satisfying $\log(\lambda(h)){\vert\chi(S)\vert}\le \log D$)? Are they all obtained by some stabilization of the elements of a finite list?
\end{question}

For example, E.\,Hironaka showed~\cite{Hironaka07} that the polynomials of smallest Mahler measure in degrees 2, 4, 6, 8, and 10 all arise as dilatations of monodromies of fibered links obtained by Hopf or trefoil plumbings on some torus links, so that the monodromies are small perturbations of some periodic surface homeomorphisms.

\begin{figure}[!ht]
	\includegraphics*[width=.4\textwidth]{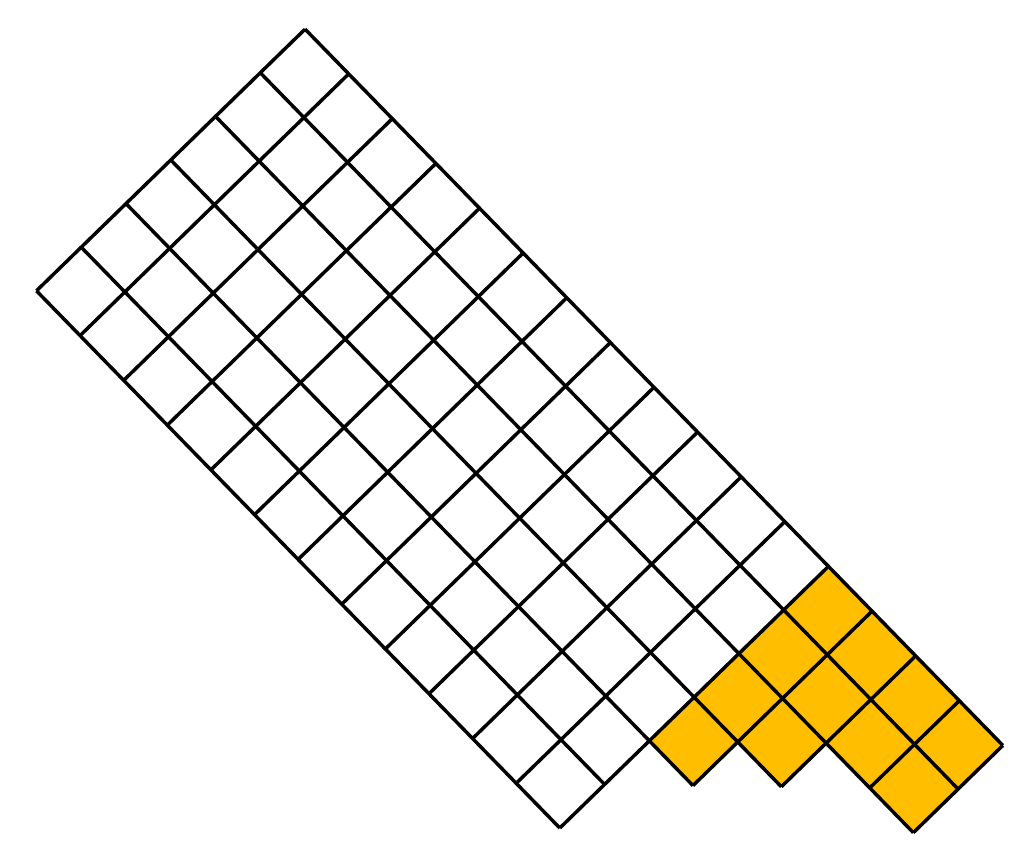}
	\caption{\small A diagram coding a knot in~$\mathcal{L}orenz_{b,k}$, for $b=6$ and $k=11$. The \emph{mixing zone} (in yellow) corresponds to those additional $k$~cells.}
  	\label{F:Allonge}
\end{figure}

What we do here is to exhibit low-dilatation families of dynamical origin, by considering certain subfamilies of the set of Lorenz knots, that is, periodic orbits of the (geometric model of the) Lorenz flow. 
These knots are fibered, so that they give rise to surface homeomorphisms, most of which are of pseudo-Anosov type.
Our statement here is a variant of Theorem~A of~\cite{DLorenz}, where we restrict our attention to subfamilies of Lorenz knots for which we obtain better bounds on the dilatation.
Denote by $\mathcal{L}orenz_{b,k}$ the set of Lorenz knots described by a hanging Young diagram (see later) made of a rectangle of width~$b$ at the bottom of which is attached a diagram with at most $k$ cells (see Figure~\ref{F:Allonge}).

\begin{theorem}
\label{T:LowDil}
The dilatation~$\lambda$ of the monodromy knot in $\mathcal{L}orenz_{b,k}$ of Euler characteristics~$\chi$ satisfies $\log(\lambda)\le \frac{b\log k}{\vert\chi\vert - k}$. 
In particular, for all $b$ and $k$, the monodromies of the elements of $\mathcal{L}orenz_{b,k}$ form a low-dilatation family.
\end{theorem}

For these families of Lorenz knots, Question~\ref{Q:QPeriodic} has a positive answer: the monodromies act like periodic homeomorphisms on a huge part of the surface (corresponding to the rectangular part of the associated Young diagram), and the non-periodicity is concentrated in a part of the surface of bounded size (corresponding to the additional $k$ cells).
Indeed, the rectangular part of the diagram corresponds exactly to a torus link, which is known to have periodic monodromy.

%





\section{Lorenz knots as iterated Murasugi sums}

Lorenz knots are defined as periodic orbits of the (geometric) Lorenz flow. 
They have been introduced and first studied by J.\,Birman and R.\,Williams~\cite{BW}. 
We refer to the original article or to~\cite{DehEns} for more details. 
Let us just mention that Lorenz knots form a family that contains all torus knots and is stable under cabling, so that it also contains all algebraic knots. 
Also, Lorenz knots are fibered, so that to each of them is canonically associated its monodromy, a homeomorphism of the genus-minimizing spanning surface. 
As Lorenz knots can be considered as perturbations of torus knots, it is natural to investigate the dilatation of the monodromies of those Lorenz knots which are hyperbolic.






\subsection{Young diagrams, Lorenz knots, and canonical spanning surfaces}
\label{S:SSurface}

\begin{figure}[!ht]
	\includegraphics*[width=.8\textwidth]{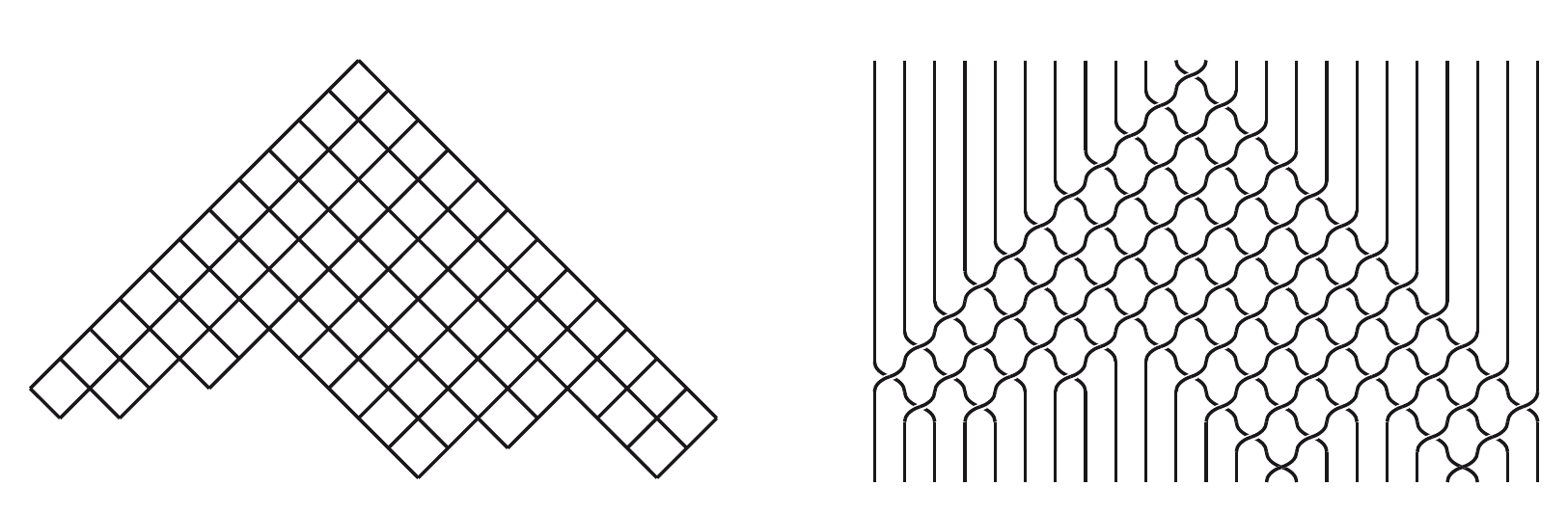}
	\caption{\small To every hanging Young diagram~$D$ (on the left), one associates a braid~$\beta(D)$ (on the right) whose closure is the Lorenz link~$K(D)$.}
  	\label{F:YoungDiagram}
\end{figure}

There are several ways of enumerating Lorenz knots and links. 
The most convenient from our point of view is using Young diagrams (introduced in this context in~\cite{DehEns}). 
The procedure is shown on Figure~\ref{F:YoungDiagram}. 

Starting from a Young diagram~$D$, one puts its bottom-left corner on top (we call this \emph{hanging position}). 
Then, by desingularizing evering intersection point into a positive braid crossing, one associates a braid, called a \emph{Lorenz braid} and denoted by $\beta(D)$.
Its closure forms a Lorenz link, that we denote by~$K(D)$. 
All Lorenz links can be obtained in this way.

Now, to the closure of every braid is associated a canonical spanning surface, obtained by gluing a disc behind every strand and a ribbon at every crossing. 
Applying this construction to~$\beta(D)$ yields a canonical spanning surface for~$K(D)$, that we denote by~$S(D)$. 
One can check that the Euler characteristics of~$S(D)$ is the number of cells of~$D$, hence denoted by~$\chi(D)$.

\begin{figure}[!ht]
	\includegraphics*[width=.3\textwidth]{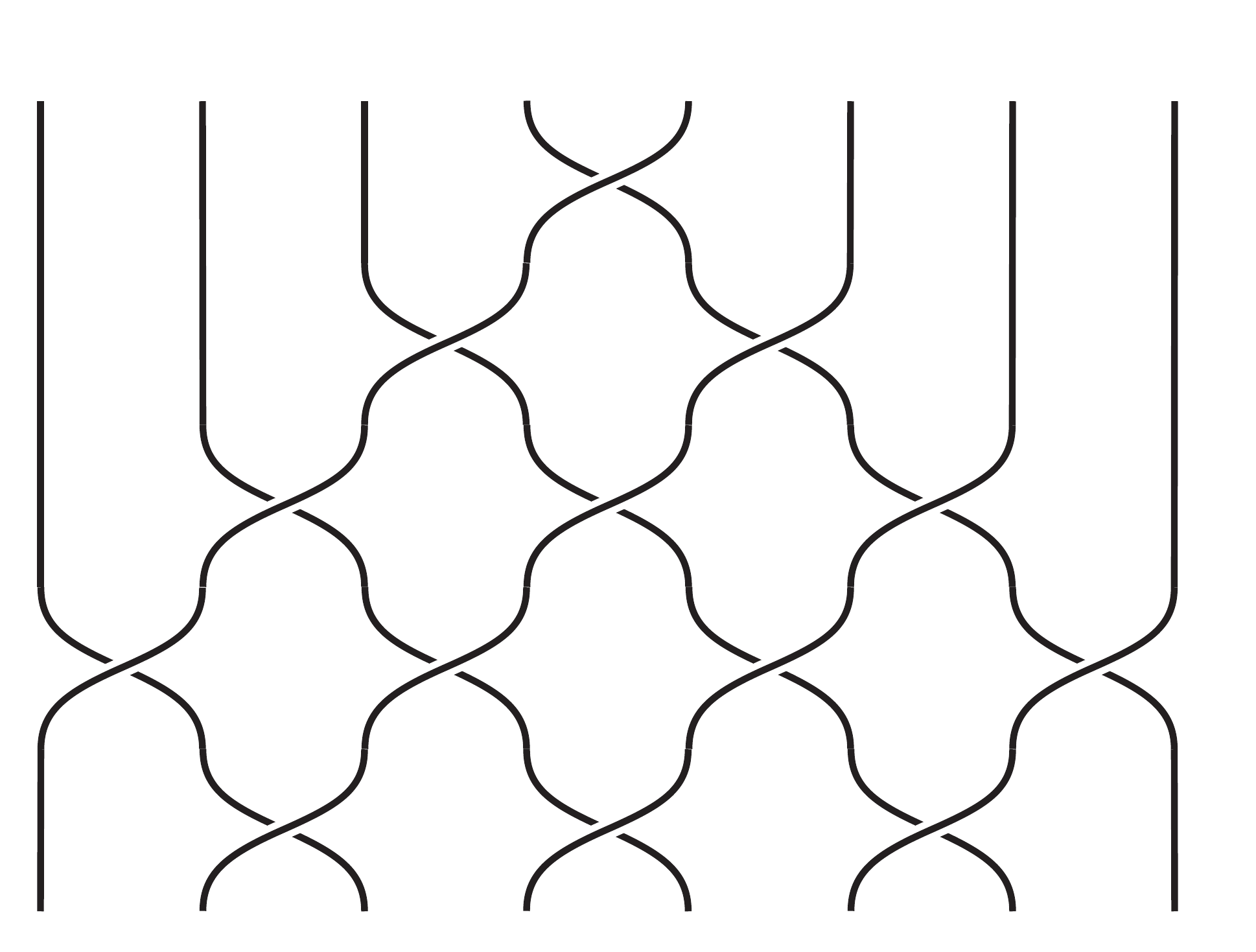}
	\quad
	\includegraphics*[width=.3\textwidth]{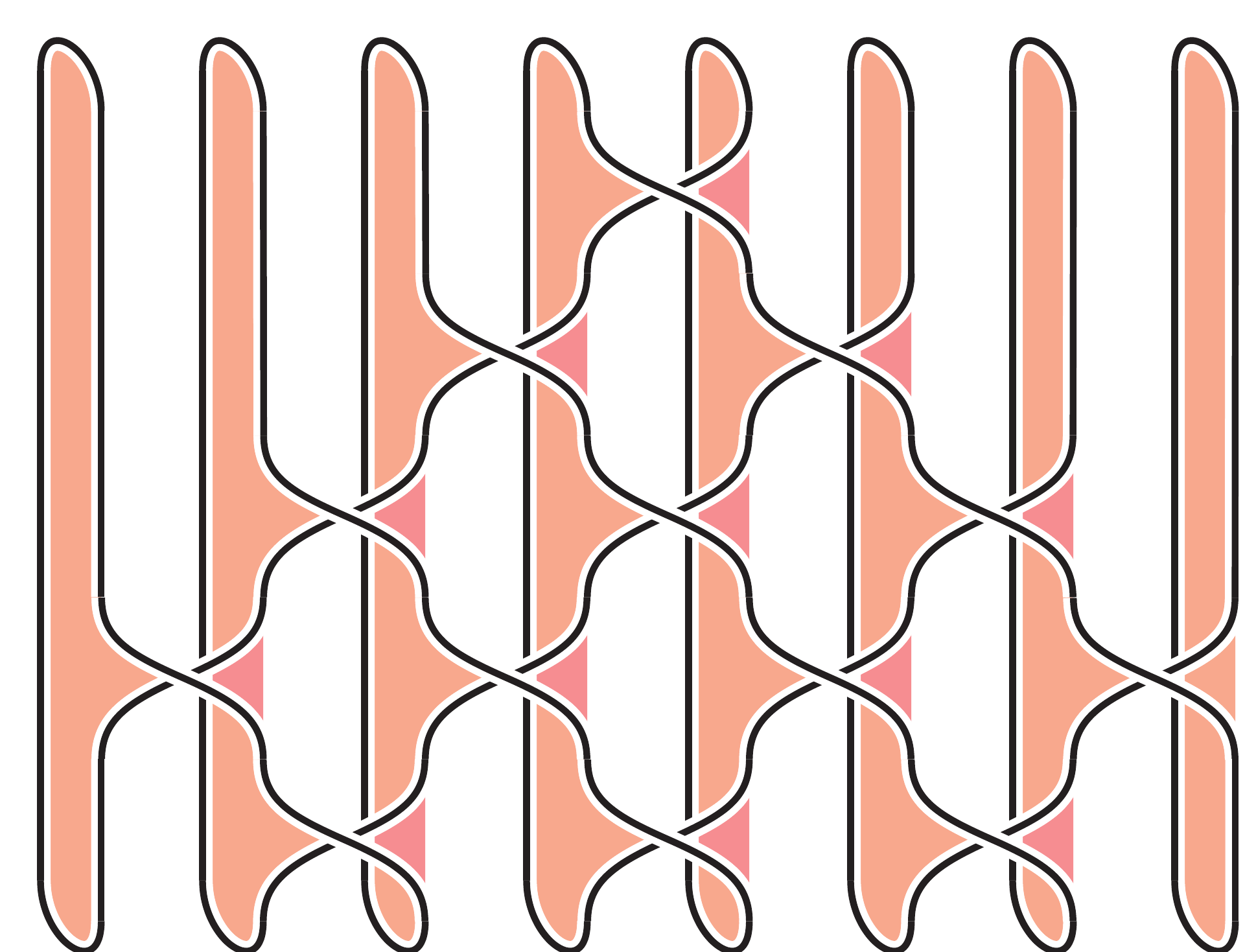}
	\quad
	\includegraphics*[width=.25\textwidth]{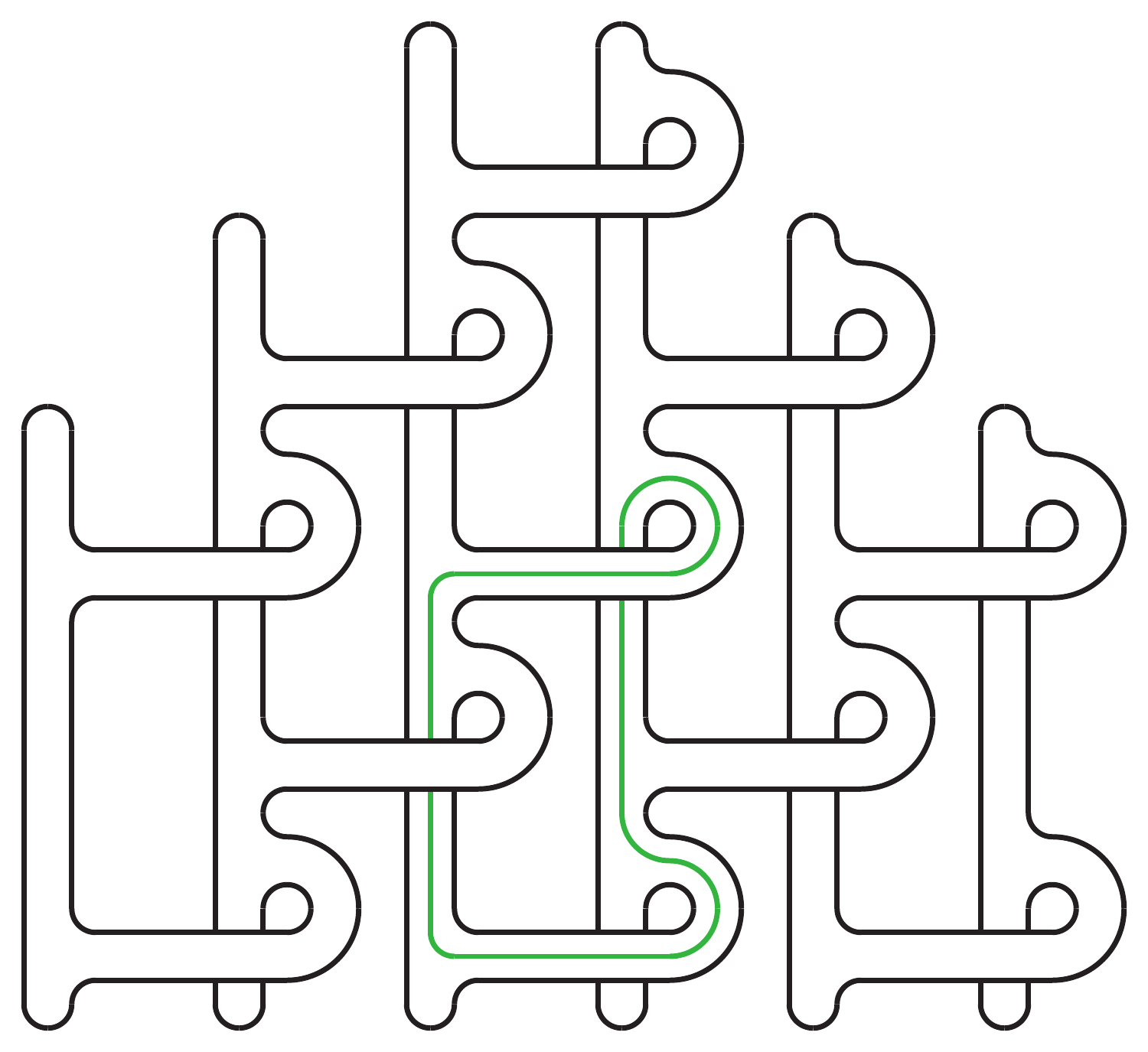}
	\caption{\small To every Lorenz braid~$\beta(D)$ (on the left), one associates a link~$K(D)$ and a canonical spanning surface~$S(D)$ (in the middle). 
	This surface can actually be immersed into the plane (following S.\,Baader~\cite{BD}, on the right). 
	In this representation, the correspondance between elementary curves (in green) and cells of the diagram is straightforward.}
  	\label{F:Surface}
\end{figure}

\subsection{Monodromy}

In this section, we describe an inductive construction of the surface $S(K)$ for every Lorenz knot~$K$, called the Murasugi sum. 
This procedure ensures that  $K$ is a fibered knot with fiber~$S(K)$, and yields a decomposition of the associated monodromy $h(K)$ as an explicit product of Dehn twists. 

By construction, for every cell~$c$ of a Young diagram~$D$ a simple close curve on~$S(D)$ that winds once around~$c$ is canonically associated. We call it a \emph{elementary curve} and denote it by~$\gamma(c)$ (see Figure~\ref{F:Surface} right).

\begin{figure}[!ht]
	\includegraphics*[width=.4\textwidth]{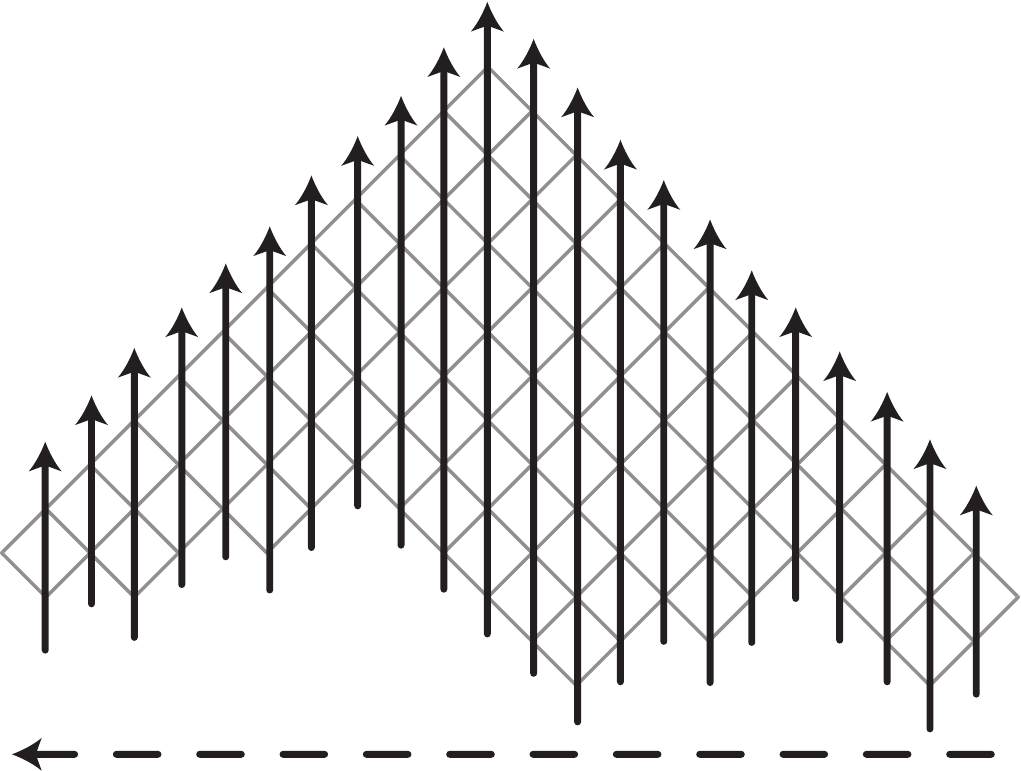}
	\caption{\small Decomposition of the monodromy $h(D)$ as the product of the Dehn twists around the $\chi(D)$ elementary curves on~$S(D)$ (that is, those curves that turns once around the cells of~$D$). The order is from right to left, and in every column from bottom to~top.}
  	\label{F:Order}
\end{figure}

\begin{proposition}
\label{P:Fiber}
Let $D$ be a Young diagram. 
Then the Lorenz link~$K(D)$ is fibered with fiber $S(D)$, and its monodromy~$h(D)$ is the product of all Dehn twists around all elementary curves of~$S(D)$, in the order prescribed on Figure~\ref{F:Order}.
\end{proposition}

Note that if $K(D)$ is a multi-component link, the fiber surface may not be unique, as well as the monodromy. 
However, if $K(D)$ is a knot, we have uniqueness of the fiber surface and of the monodromy homeomorphism.

We will only sketch the proof of Proposition~\ref{P:Fiber} and refer to the survey~\cite{DehEns} for more details. 
The starting point is the 2-component Hopf link, which is the Lorenz link associated to the Young diagram with one cell only. 
The 2-component Hopf link is known to be fibered, the fiber surface being a twisted annulus that we call the \emph{Hopf annulus}, and the monodromy being a right-handed Dehn twist.

\begin{figure}[!ht]
	\includegraphics*[width=.15\textwidth]{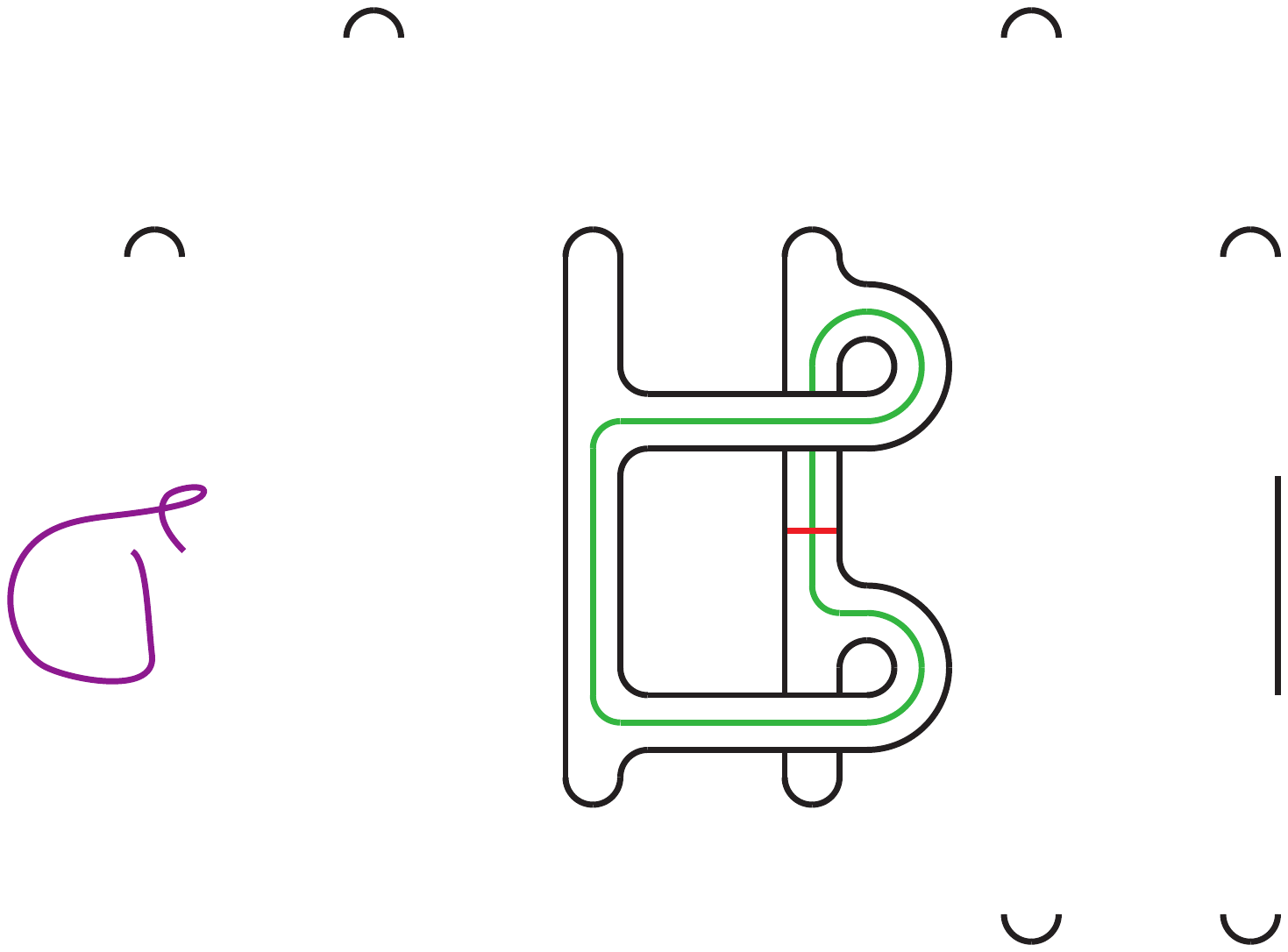}
	\qquad
	\includegraphics*[width=.2\textwidth]{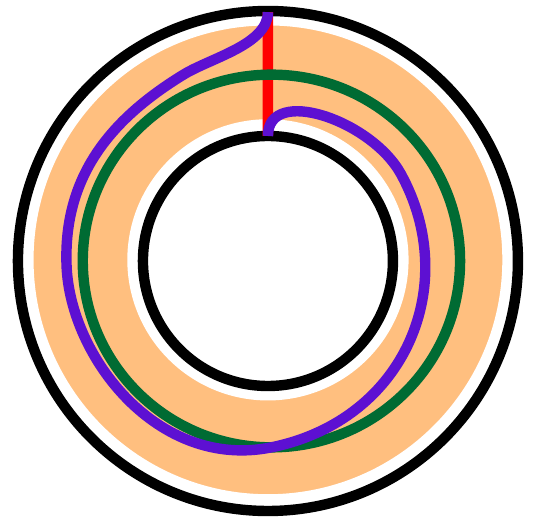}
	\caption{\small A Hopf annulus in~$\Sph^3$(on the left), with the associated elementary curve (in green). The action of the monodromy on the annulus (on the right, seen on an abstract annulus) is a Dehn twist on the green curve: the red segment is sent to the purple one.}
  	\label{F:Hopf}
\end{figure}

\begin{figure}[!ht]
	\includegraphics*[width=.35\textwidth]{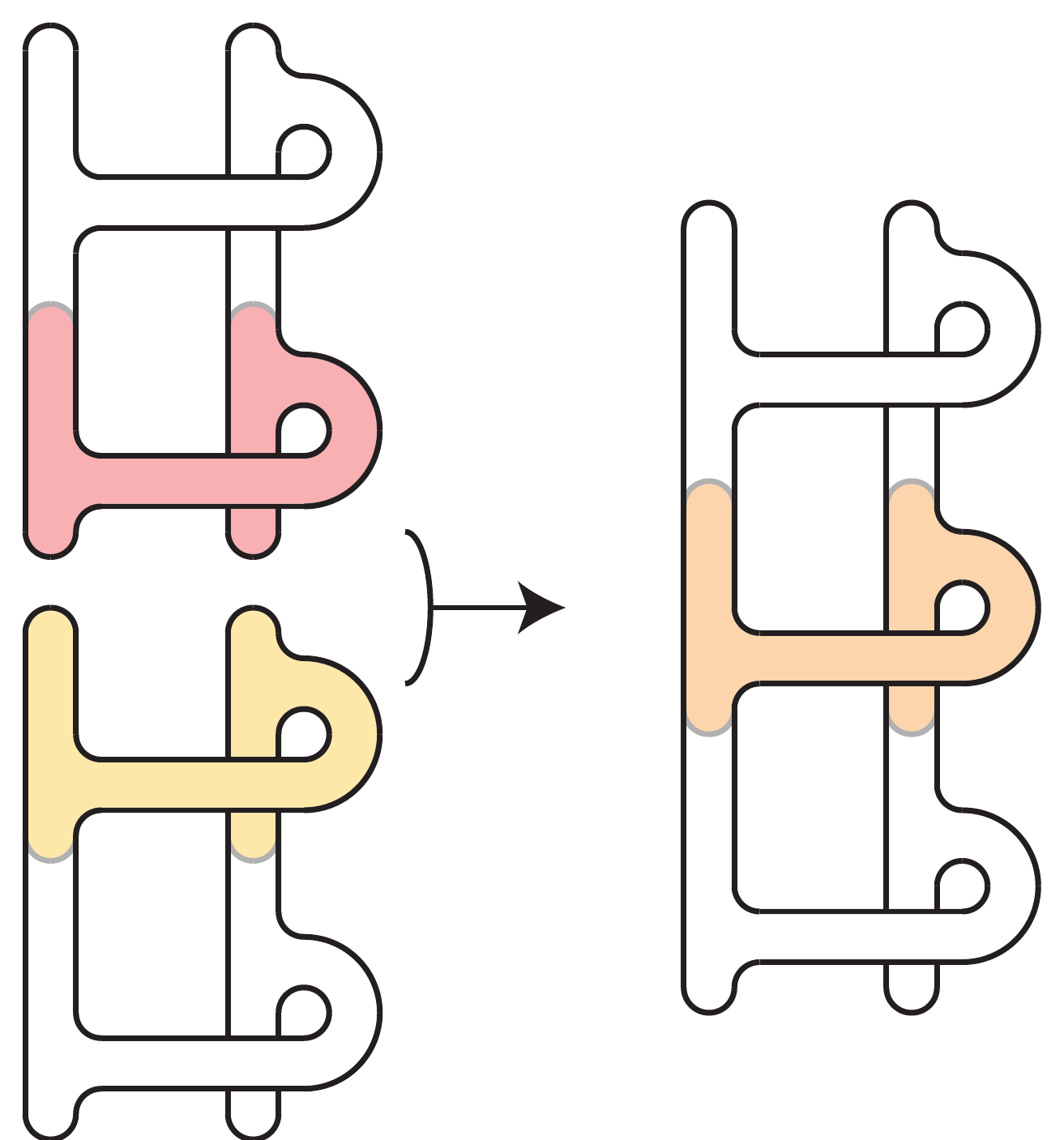}

	\includegraphics*[width=.4\textwidth, angle = 90]{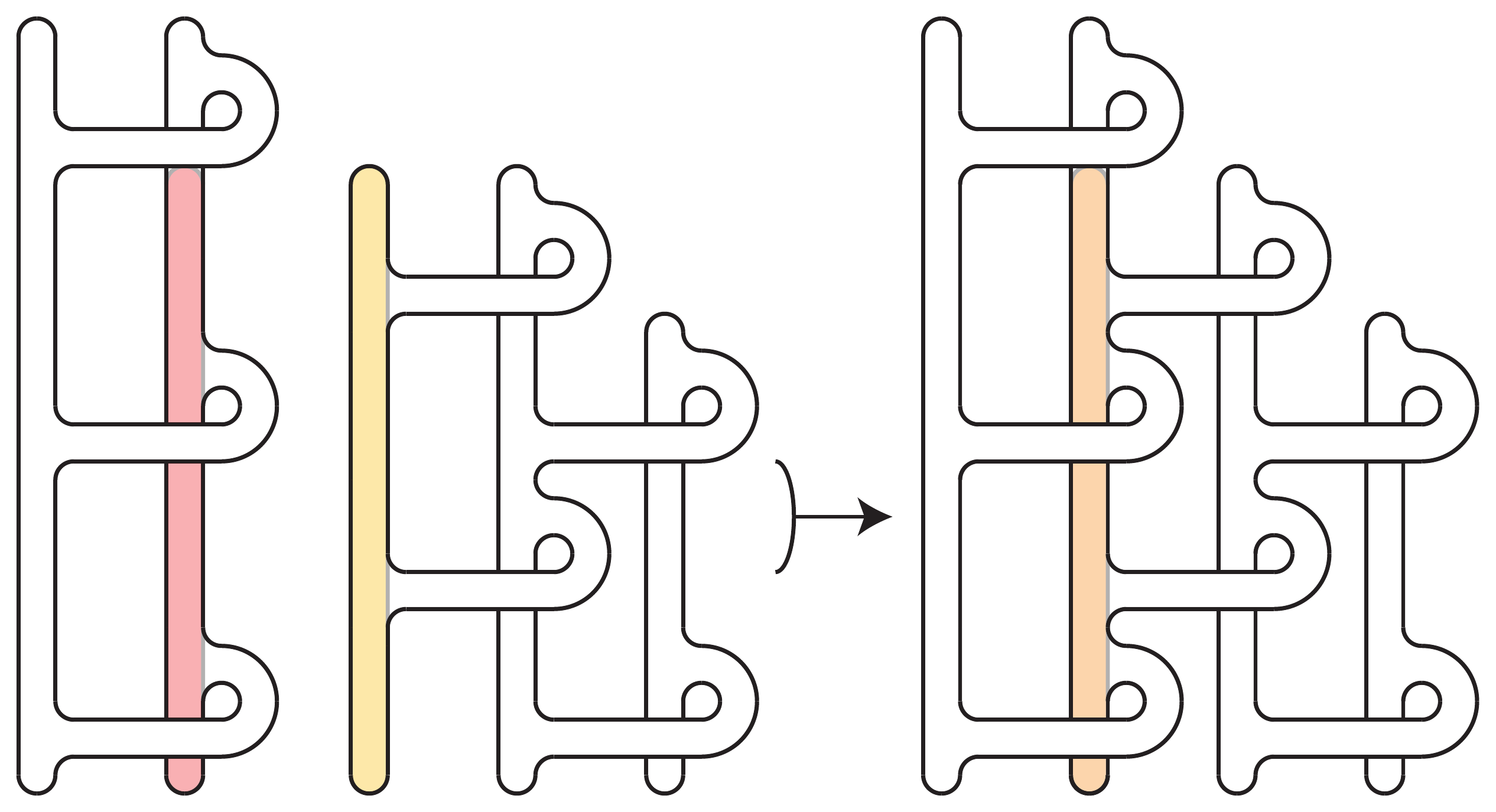}
	\caption{\small Two examples of Murasugi sums of canonical surfaces associated to positive braids. Observe that in the second example the left summand is glued on top.}
  	\label{F:Murasomme}
\end{figure}

The induction step for proving Proposition~\ref{P:Fiber} is done using the \emph{Murasugi sum} of surfaces~\cite{Murasugi}. 
This is an operation that takes two surfaces with boundary $S_1, S_2$, depends on a choice of a $2n$-gon in each of them, and associates a new surface with boundary~$S_1\sharp S_2$ that contains $S_1$ and $S_2$ as subsurfaces (see Figure~\ref{F:Murasomme}).
This operation preserves the fibered character, in the following sense: 
if $S_1, S_2$ are two fibered surfaces in~$\Sph^3$ with monodromies~$h_1, h_2$, then the Murasugi sum~$S_1\sharp S_2$, where $S_1$ is glued on top, is fibered with monodromy $h_1\circ h_2$ (see the proof of D.\,Gabai~\cite{Gabai} or an expanded version in~\cite{DehEns}).

In particular, Murasugi gluing a Hopf annulus to a fibered surface yields another fibered surface.
In this way, starting from the canonical Seifert surface associated to a hanging Young diagram~$D$, we obtain that the surface associated to the diagram $D'$ obtained from~$D$ by adding a cell on the bottom-right border of~$D$ is also fibered. 
Moreover, the monodromy associated to a Hopf annulus is a right-handed Dehn twist so that the monodromy associated to the surface~$S(D)$ is a product of Dehn twists along the cores of the glued Hopf annuli. 
The order of the product is determined by the order of the gluing.
The latter needs to preserve the respective positions of the Hopf annuli, namely one should glue first an annulus that is on top of another one.
The order given on Figure~\ref{F:Order} obeys this constaint.
This completes the (sketch of) proof of Proposition~\ref{P:Fiber}.

\subsection{Action of the monodromy on elementary curves}

The dilatation of a pseudo-Anosov homeomorphism can be read on its action on curves. 
So, in order to bound the dilatation, one should bound the stretching of curves under the homeomorphism.
Cutting the canonical surface~$S(D)$ associated to a diagram~$D$ along all elementary curves reduces~$S(D)$ to a neighborhood of its boundary, so that elementary curves contain all the information on~$S(D)$. 
In particular for Lorenz knots, it is enough to estimate the stretching of elementary curves under~$h(D)$ is order to control the dilatation of~$h(D)$.

Now come the two key observations. 
For some orientation reason, the second observation works only when considering $h^{-1}(D)$ instead of~$h(D)$.
Therefore we consider the inverse of the monodromy, which makes little difference. 

\begin{figure}[!ht]
	\includegraphics*[width=.35\textwidth]{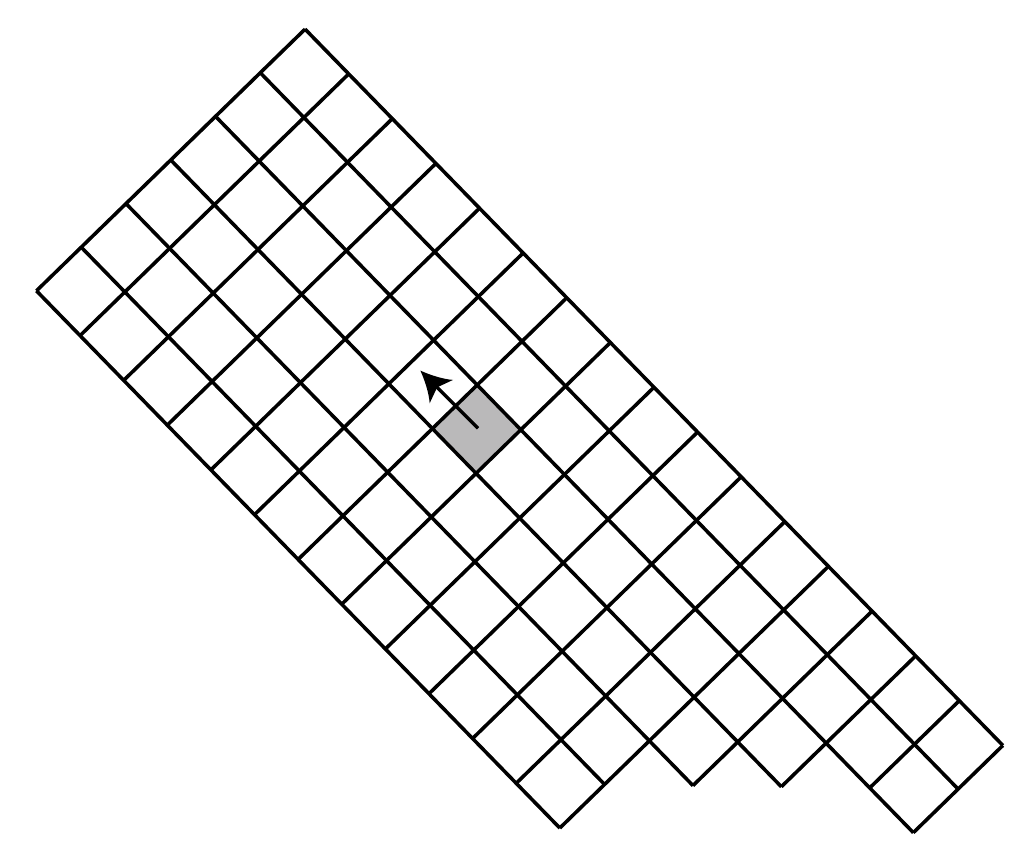}
	\caption{\small An internal cell and its image under the (inverse of the) monodromy.}
  	\label{F:Intern}
\end{figure}

We say that a cell~$c$ of a hanging Young diagram is \emph{internal} if there is a cell, say~$c'$, in North-West position with respect to~$c$ (see Figure~\ref{F:Intern}). Otherwise it is called~\emph{external}.

\begin{lemma}
\label{L:ImInt}
Assume that $D$ is a Young diagram, $L(D)$ is the associated Lorenz link, $S(D)$ is the canonical Seifert surface for~$L(D)$, and $h(D)$ is the associated monodromy.
Let $d$ be an internal cell of~$D$ and $a$ be the cell in NW position with respect to~$c$. 
Then we have $h(D)^{-1}(\gamma(d))=\gamma(a)$.
\end{lemma}

\begin{figure}[!ht]
	\includegraphics*[width=\textwidth]{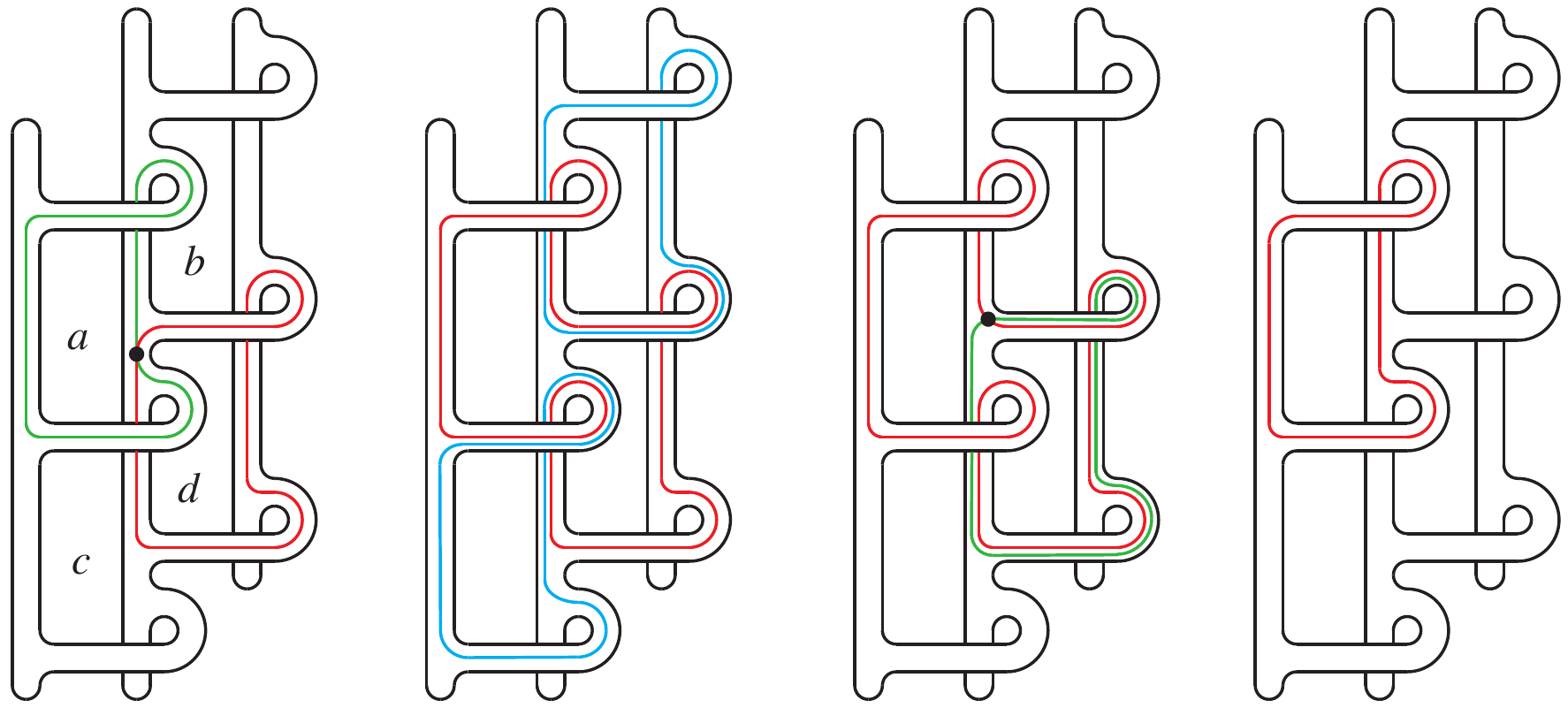}
	\caption{\small Proof of Lemma~\ref{L:ImInt}. The image of the elementary curve~$\gamma(d)$ associated to an internal cell~$d$ under the inverse of the monodromy $h(D)^{-1}$: first the Dehn twists associated to cells that are distant from~$d$ do not modify~$\gamma(d)$. Then it is changed by the Dehn twist around $a$ into a curve that encircles both~$a$ and $d$. This curve is then unchanged (in particular it does not intersect the two blue curves~$\gamma(b)$ and $\gamma(c)$ on the second picture). Finally it is changed by the Dehn twist around $d$ into~$\gamma(a)$. Subsequent twists to not modify it any more.}
  	\label{F:MonInt}
\end{figure}

The proof is displayed on Figure~\ref{F:MonInt}, where the successive images of the curve~$\gamma(d)$ under consecutive Dehn twists are depicted (see also~\cite{BD}). 

In order to fully control~$h(D)^{-1}$, we need to know what happens to external cells when iterating (backwards) the monodromy. 
For a general Lorenz link, the behaviour is rather hard to control (this is the reason of the heavy computations in~\cite{DLorenz}). 
However, if we suppose that the diagram we are considering lies in~$\mathcal{L}orenz_{b,k}$, things become simpler.
In particular the image of an elementary curve corresponding to an external cell is not so simple, but its second image is. 

For $D$ a diagram in~$\mathcal{L}orenz_{b,k}$, we call \emph{mixing zone} of~$D$ the set of those $k$ cells that are outside the main rectangle of~$D$ (see Figure~\ref{F:Allonge}). 
We also assume that we have an auxiliary metric on~$S$ for which all elementary curves have length at most~$1$.

\begin{lemma}
\label{L:ImExt}
Assume that $D$ is a Young diagram in~$\mathcal{L}orenz_{b,k}$, and that $h(D)$ is the monodromy associated to the canonical surface~$S(D)$.
Let $c$ be an external cell of~$D$.
Then $h(D)^{-2}(\gamma(c))$ is a curve of length at most $k$ that lies entirely in the mixing zone.
\end{lemma}

\begin{figure}[!ht]
	\includegraphics*[width=.225\textwidth, angle = 270]{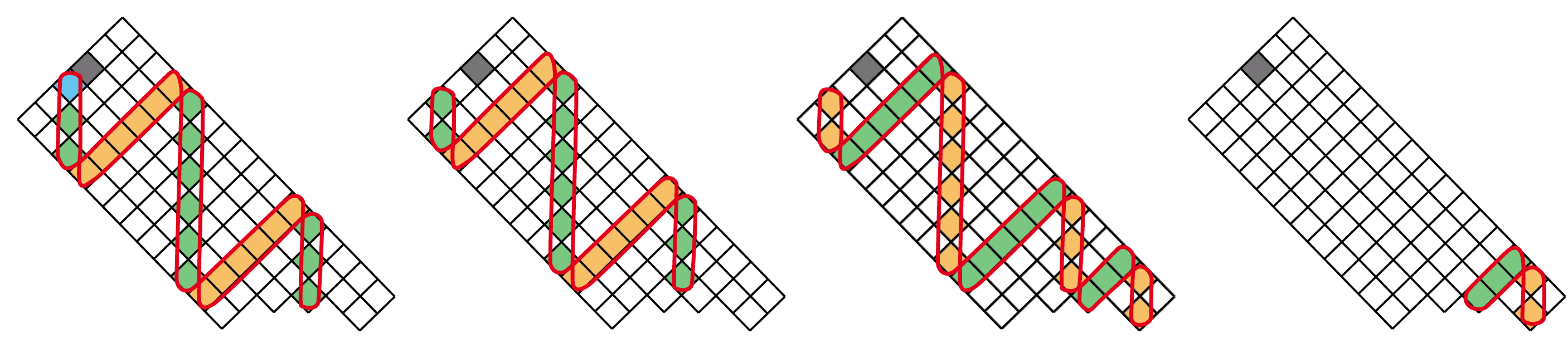}
	\caption{\small Proof of Lemma~\ref{L:ImExt}: an external cell~$c$ and its image under $h(D)^{-2}$. On the left, $h(D)^{-1}(\gamma(c))$ is a curve that turns positively around blocks of orange cells and negatively around blocks of green (and blue) cells. On the center left, the image under $h(D)^{-1}$ of the curve $h(D)^{-1}(\gamma(c))$, except the part that winds around the blue cell. On the center right, the image under $h(D)^{-1}$ of the part of the curve $h(D)^{-1}(\gamma(c))$ that winds around the blue cell. On the right, the concatenation of those two parts is the curve~$h(D)^{-2}(\gamma(c))$, it is a curve that only winds around some cells of the mixing zone of~$D$.}
  	\label{F:Extern}
\end{figure}

The proof is depicted on Figure~\ref{F:Extern}. 
The idea is that, with arguments similar to the proof of Lemma~\ref{L:ImInt}, one can describe the curve $h(D)^{-1}(\gamma(c))$: it is the concatenation of one external curve, and many internal curves. When iterating $h(D)^{-1}$ once more, the different contributions cancel, except in the mixing zone.

\subsection{Proof of Theorem~\ref{T:LowDil}}
Assume that~$D$ is a Young diagram in~$\mathcal{L}orenz_{b,k}$, that~$S(D)$ is the associated canonical surface, and that $h(D)$ is the corresponding monodromy. 
Denote by~$l(D)$ the length of the long rectangle in~$D$ (that is, the complement of the mixing zone). 
We also take an auxiliary metric on~$S(D)$ for which all elementary curves have length at most~1.

Let $c$ be an arbitrary cell in the mixing zone on~$S(D)$. 
By Lemma~\ref{L:ImInt}, the $l(D)$ first images of the curve~$\gamma(c)$ under~$h(D)^{-1}$ all correspond to internal cells, hence have length one. 
After a few more iteration, the image is then an elementary external curve, and after two more iterations, it is a curve in the mixing zone of length at most~$k$. 
Then the process goes on: the next $l$ iterations yield a curve of length at most~$k$.
Summarizing, the length of $h(D)^{-n}(\gamma(c))$ grows by a factor at most~$k$ every $l(D)$ steps. 
Therefore, the growth rate of~$\gamma(c)$ is bounded by~${\log k}/{l(D)}$.

Now, the same argument works for any cell, not just in the mixing zone, except that the initial dilatation arises earlier. 
But this does not change the growth rate, hence $h(D)^{-1}$ asymptotically streches all curves on~$S(D)$  by a factor at most~${\log k}/{l(D)}$.

Finally, an elementary computation shows that the Euler characteristics of~$S(D)$ is~$b\cdot l(D)+k$, therefore the dilatation is smaller than ${b\log k}/({\vert\chi(D)\vert-k})$.


\bibliographystyle{alpha}

\end{document}